\documentclass[11pt]{amsart}
\usepackage{amssymb}

\newcommand{\script}[1]{\mathcal{#1}}
\newcommand{\Lie}{\mathcal}
\newcommand{\Gs}[1]{{#1}_S^{\scriptscriptstyle(1)}}
\newcommand{\GS}{\Gs{\Ga}}
\newcommand{\bigtimes}{\mathop{\hbox{\Large$\times$}}}
\newcommand{\gi}{G_\infty}
\newcommand{\gio}{G_\infty^\circ}
\newcommand{\finite}{{\textnormal{f}}}
\newcommand{\gf}{G_\finite}
\newcommand{\G}{\closure{\Ad G}}
\newcommand \Gi[1]{\closure{\Ad G_{#1}}}
\newcommand \Go{\closure{\Ad G^{\circ}}}
\newcommand{\Z}[1]{\closure{\Ad_{G}#1}}
\newcommand \Zi[2]{\closure{\Ad_{G_#1} #2}}
\newcommand{\theints}{{\script O(S)}}
\newcommand{\Sints}{{\integer(S)}}
\newcommand{\Ga}{\mathbb G}
\newcommand{\Ta}{\mathbb T}
\newcommand{\Ua}{\mathbb U}
\newcommand{\La}{\mathbb L}
\newcommand{\therad}{R}
\newcommand{\GoS}[1]{{\closure{{#1}_{\theints}}}}
\newcommand{\GOS}{\GoS{\Ga}}
\newcommand{\GL}{\operatorname{GL}}
\newcommand{\SO}{\operatorname{SO}}
\newcommand{\Ad}{\operatorname{Ad}}
\newcommand{\Rad}{\operatorname{Rad}}
\newcommand{\Aut}{\operatorname{Aut}}
\newcommand{\nil}{\operatorname{nil}}
\newcommand{\closure}{\overline}
\newcommand{\cover}{\widetilde}
\newcommand{\graph}{\operatorname{graph}}
\newcommand{\iso}{\cong}
\newcommand{\semiprod}{\ltimes}
\newcommand{\prodsemi}{\rtimes}
\newcommand{\real}{\mathord{\mathbb R}}
\newcommand{\rational}{\mathord{\mathbb Q}}
\newcommand{\complex}{\mathord{\mathbb C}}
\newcommand{\integer}{\mathord{\mathbb Z}}
\newcommand{\torus}{\mathord{\mathbb T}}

\newcommand \stackarrow[1]{\setbox0\hbox{ \ $\scriptstyle#1$ \ }
 \setbox1\hbox{$\rightarrow$}\setbox2\hbox to
   \wd0{\rightarrowfill}\ht2=\ht1
 \mathrel{\mathop{\copy2}\limits^{\copy0}}}

\newcommand \CASE[2]{\medskip \noindent \emph{Case \textnormal{#1.} #2}
   \ignorespaces}
\newcommand \step[2]{\medskip \noindent \emph{Step \textnormal{#1.} #2}
   \ignorespaces}

\renewcommand{\see}[1]{\textnormal{(}see~\textnormal{\protect\ref{#1})}}
\newcommand{\seemore}[2]{\textnormal{(}see#1~\textnormal{\protect\ref{#2})}}
\newcommand{\pref}[1]{(\protect\ref{#1})}

\renewcommand{\thmhead}[3]{%
 \thmname{#1}\thmnumber{ #2}\thmnote{ {\the\theoremnotefont#3}}}

\renewcommand{\swappedhead}[3]{%
 \textnormal{(\thmnumber{#2})}\thmname{ #1}\thmnote{
 {\the\theoremnotefont#3}}}
\swapnumbers

\newtheorem{thm}{Theorem}[section]
\newtheorem{mainthm}[thm]{Main Theorem}
\newtheorem{lem}[thm]{Lemma}
\newtheorem{cor}[thm]{Corollary}
\newtheorem{prop}[thm]{Proposition}

\theoremstyle{definition}

\newtheorem{defn}[thm]{Definition}
\newtheorem{eg}[thm]{Example}
\newtheorem{assump}[thm]{Assumption}

\theoremstyle{remark}

\newtheorem{rem}[thm]{Remark}        
\newtheorem{ack}[thm]{Acknowledgment}  

\newbox\sectionS \setbox\sectionS\hbox{\S}

\newcommand{\refjour}{} %% make sure the names aren't already used
\newcommand{\refbook}{}

\long\def\refjour[#1]#2 #3: #4. #5.. #6 (#7) #8 \par
 {\bibitem{#2} #3, #4, \emph{#5} {\bf #6} (#7), #8.\par}
\long\def\refbook[#1]#2 #3: #4. #5, #6, #7 \par
 {\bibitem{#2} #3, ``#4,'' #5, #6, #7.\par}

\begin{document}

\title[Superrigidity of $S$-arithmetic groups]
{Archimedean superrigidity of \\ solvable $S$-arithmetic groups}

\author{Dave Witte}

\address{Department of Mathematics, Oklahoma State University,
 Stillwater, OK 74078}

\email{dwitte@math.okstate.edu}

\date{\today}

%%%%%%%%%%%%%%%%%%%%%%%%%%%%%%%%%%%%%%%%%%%%%%%%%%%%%%%%%%%%
%%% Publisher: you'll want to delete the following line. %%%
%%%%%%%%%%%%%%%%%%%%%%%%%%%%%%%%%%%%%%%%%%%%%%%%%%%%%%%%%%%%
\thanks{To appear in \emph{Journal of Algebra}. (Submitted March 1996.)}

\begin{abstract}
 Let $\Ga$ be a connected, solvable linear algebraic group over a
number field~$K$, let $S$ be a finite set of places of~$K$ that
contains all the infinite places, and let $\theints$ be the ring of
$S$-integers of~$K$. We define a certain closed subgroup~$\GOS$ of
$\Ga_S = \prod_{v \in S} \Ga_{K_v}$ that contains $\Ga_{\theints}$,
and prove that $\Ga_{\theints}$ is a superrigid lattice in~$\GOS$, by
which we mean that finite-dimensional representations $\alpha\colon
\Ga_{\theints} \to \GL_n(\real)$ more-or-less extend to
representations of~$\GOS$.

 The subgroup~$\GOS$ may be a proper subgroup of~$\Ga_S$ for only two
reasons. First, it is well known that $\Ga_{\theints}$ is not a
lattice in~$\Ga_S$ if $\Ga$ has nontrivial $K$-characters, so one
passes to a certain subgroup $\GS$. Second, $\Ga_{\theints}$ may fail
to be Zariski dense in $\GS$ in an appropriate sense; in this sense,
the subgroup $\GOS$ is the Zariski closure of~$\Ga_{\theints}$
in~$\GS$.

Furthermore, we note that a superrigidity theorem for many
non-solvable $S$-arithmetic groups can be proved by combining our main
theorem with the Margulis Superrigidity Theorem.
 \end{abstract}

\maketitle

\section{Introduction} \label{introduction}

Let $\Ga$ be a solvable linear algebraic group defined
over~$\rational$.  The author recently proved that if the arithmetic
subgroup~$\Ga_{\integer}$ is Zariski dense, then it is a superrigid
lattice in~$\Ga_{\real}$ \see{my-thm}, in the sense that any
finite-dimensional representation $\alpha\colon \Ga_{\integer} \to
\GL_n(\real)$ more-or-less extends to a representation
of~$\Ga_{\real}$. (A precise definition of superrigidity appears in
Defn.~\ref{super-def} below.) We now prove an appropriate
generalization of this result for $S$-arithmetic subgroups, in place
of arithmetic subgroups (see~\ref{super-Q} and~\ref{main-thm}).

\begin{defn}[{(cf.~\cite[Thm.~2, p.~2]{MargBook})}] \label{super-def}
 Let $\Lambda$ be a subgroup of a topological group~$G$. We say that $\Lambda$
is \emph{superrigid} in~$G$ if, for every continuous homomorphism $\alpha
\colon \Lambda \to \GL_{n}(\real)$,
 there are
 \begin{itemize}
 \item a finite-index open subgroup~$\Lambda_1$ of~$\Lambda$,
 \item a finite-index open subgroup~$G_{1}$ of~$G$, containing~$\Lambda_1$, and
 \item a finite, normal subgroup~$F$ of
$\closure{\Lambda_1^\alpha}$, where $\closure{\Lambda_1^\alpha}$ is
the (almost-)Zariski closure of $\Lambda_1^\alpha$ in
$\GL_n(\real)$,
 \end{itemize}
 such that the induced homomorphism $\alpha_F\colon \Lambda_1 \to
\closure{\Lambda_1^\alpha}/F$ extends to a continuous homomorphism
 $\sigma \colon G_1 \to \closure{\Lambda_1^{\alpha}}/F$.

 Note that the only superrigidity under consideration is
\emph{archimedean} superrigidity. That is, all representations are
over~$\real$ (or~$\complex$). Our results say nothing about
representations over other local fields.
  \end{defn}

\begin{rem}
 By considering restrictions and induced representations, one may
show, for any finite-index subgroup~$\Lambda'$ of~$\Lambda$, that
$\Lambda$~is superrigid in~$G$ if and only if $\Lambda'$ is
superrigid in~$G$.
 \end{rem}

\begin{assump}
For simplicity, we assume that all algebraic groups in this paper are
linear. (More precisely, they are subgroups of some special linear
group.) Thus, if $\Ga$ is an algebraic group over~$\rational$, then
$\Ga_{\integer}$ is well defined, and is a subgroup
of~$\Ga_{\rational}$. Without such an assumption, the definition of
arithmetic subgroups of~$\Ga$ is slightly more involved \cite[\S7.11,
p.~49]{Borel-arith}.
 \end{assump}

\begin{thm}[{(\cite[Thm.~6.8]{Witte-super}, see
also~\S\ref{easy-pf})}] \label{my-thm}
 Let\/ $\Ga$ be a solvable algebraic group over\/~$\rational$. If\/
$\Ga_{\integer}$ is Zariski dense in\/~$\Ga$, then\/\/ $\Ga_{\integer}$ is a
superrigid lattice in\/~$\Ga_{\real}$.
 \end{thm}

If $\Ga$ has no nontrivial characters defined over~$\rational$ (and is
connected), then this theorem has a natural generalization to
$S$-arithmetic subgroups.

\begin{defn}
 Let $\Ga$ be a connected algebraic group over~$\rational$, and let $S$ be a
finite set of places of~$\rational$, containing the infinite place.
 Define
 \begin{itemize}
 \item $\Sints = \{ x \in \rational \mid \hbox{$||x ||_v \le 1$ for all
places $v \not\in S$} \}$;
 \item $\rational_v =$ the completion of~$\rational$ at the place~$v$; and
 \item $\Ga_S = \bigtimes_{v \in S} \Ga_{\rational_v} $.
 \end{itemize}
 Via the natural diagonal embedding of $\Ga_{\rational}$ in~$\Ga_S$, we may
think of~$\Ga_{\Sints}$ as a subgroup of~$\Ga_S$.
 \end{defn}

\begin{thm} \label{super-Q}
 Let\/ $\Ga$ be a connected, solvable algebraic group over\/~$\rational$,
and let $S$ be a finite set of places of\/~$\rational$, containing the
infinite place. If\/ $\Ga$ has no nontrivial characters defined
over\/~$\rational$, and\/ $\Ga_\Sints$ is Zariski dense in\/~$\Ga$, then\/
$\Ga_\Sints$ is a superrigid lattice in\/~$\Ga_S$. \qed
 \end{thm}

If $\Ga$ does have nontrivial $\rational$-characters, then $\Ga_\Sints$ is
not a lattice in~$\Ga_S$. However, it is well known that $\Ga_\Sints$ is a
lattice in a certain subgroup~$\GS$ \see{GOS-def}
and~\ref{is-lattice}), and our main theorem implies that $\Ga_\Sints$ is
superrigid in~$\GS$. More generally, our main theorem shows that an
analogous result holds for solvable algebraic groups over any number
field, but, in this generality, it may be necessary to replace~$\GS$ with
a smaller group~$\GOS$, which we now define \see{S-defs}. The idea is that
$\GOS$ should be the almost-Zariski closure \see{almost-def} of
$\Ga_{\theints}$ in~$\GS$, but the situation is complicated by the fact
that $\GS$ is not an algebraic group.

\begin{defn}[{\cite[Defn.~3.2]{Witte-super}}]
 A subgroup~$A$ of $\GL_n(\real )$ is said to be \emph{almost Zariski closed} if
there is a Zariski closed subgroup~$B$ of $\GL_n(\real )$, such that $B^\circ
\subset A \subset B$, where $B^\circ$ is the identity component of~$B$ in the
topology of $\GL_n(\real)$ as a $C^\infty$~manifold (not the Zariski
topology). There is little difference between being Zariski closed and almost
Zariski closed, because $B^\circ$ always has finite index in~$B$.
 \end{defn}

\begin{defn}[{\cite[Defn.~3.6]{Witte-super}}] \label{almost-def}
 Let $A$ be a subgroup of $\GL_n(\real)$. The \emph{almost-Zariski
closure}~$\closure{A}$ of~$A$
 is the unique smallest almost-Zariski closed subgroup that contains~$A$.
 In particular, if $A$ is a subgroup of a Lie group~$G$, we use $\Z{A}$ to
denote the almost-Zariski closure of $\Ad_G A$ in the real
algebraic group $\Aut(\Lie G)$, where $\Lie G$ is the Lie algebra
of~$G$.
 \end{defn}

\begin{defn} \label{S-defs}
 Let $\Ga$ be a connected algebraic group over a number field~$K$, and let $S$
be a finite set of places of~$K$, containing the set~$S_\infty$ of infinite
places.
 Define
 \begin{itemize}
 \item $\theints = \{ x \in K \mid \hbox{$||x ||_v \le 1$ for all
places $v \not\in S$} \}$;
 \item $K_v =$ the completion of~$K$ at the place~$v$; 
 \item $\Ga_\infty = \bigtimes_{v \in S_\infty} \Ga_{K_v} $;
 \item $\Ga_{\textnormal{f}} = \bigtimes_{v \in S - S_\infty}
\Ga_{K_v}$;
 \item $\Ga_S = \Ga_\infty \times \Ga_{\textnormal{f}}$;
 \item $\GS = 
 \{ x \in \Ga_S \mid
 \hbox{$\prod_{v \in S} || \chi(x_v) ||_v = 1$, for all $K$-characters~$\chi$
of~$\Ga$} \}$
 \cite[pp.~263--264]{Plat-Rap}; and
 \item $\GOS = \{ x \in \GS \mid \Ad_{\Ga_\infty} x \in
\closure{\Ad_{\Ga_\infty} \Ga_\theints} \} $. \label{GOS-def}
 \end{itemize}
 See Ex.~\ref{GOS-proper} for an example where $\GOS$ is a proper
subgroup of $\GS$.
 \end{defn}

\begin{mainthm}[{(see \S\ref{main-pf-sect})}] \label{main-thm}
 Let\/ $\Ga$ be a connected, solvable algebraic group over a number
field~$K$, and let $S$ be a finite set of places of~$K$, containing all the
infinite places. If\/ $\Ga_\theints$ is Zariski dense in\/~$\Ga$, then\/
$\Ga_\theints$ is a superrigid lattice in\/ $\GOS$.
 \end{mainthm}

The following corollary states that the question of whether $S$-arithmetic
subgroups of a given algebraic group~$\Ga$ are superrigid can essentially be
reduced to the same question about the $S$-arithmetic subgroups of the maximal
semisimple quotient of~$\Ga$. The point is that the main theorem implies that
solvable $S$-arithmetic groups are superrigid, which means the radical
of~$\Ga$ is under control---all that remains is the semisimple part of~$\Ga$.

\begin{cor}[{(see \S\ref{non-solv-sect})}] \label{main-cor}
 Let\/ $\Ga$ be a connected algebraic group over a number field~$K$, and let
$S$ be a finite set of places of~$K$, containing the infinite places. If\/
$\Ga_\theints$ is Zariski dense in\/~$\Ga$,
 and the image of\/ $\Ga_\theints$ in\/ $\GoS{(\Ga/\Rad\Ga)}$ is
superrigid, then\/ $\Ga_\theints$ is a superrigid lattice in\/~$\GOS$.
 \end{cor}

G.~A.~Margulis has shown that $S$-arithmetic subgroups of semisimple
groups of higher $S$-rank are superrigid.

\begin{thm}[{(Margulis \cite[B(iii), p.~259]{MargBook})}]
 Let\/ $\Ga$ be a connected, semisimple algebraic group over a number field~$K$,
and let $S$ be a finite set of places of~$K$, containing the infinite
places. Assume the $S$-rank of every simple factor of\/~$\Ga$ is at least
two. Then\/ $\Ga_\theints$ is a superrigid lattice in\/~$\Ga_S$. \qed
 \end{thm}

\begin{cor}
 Let\/ $\Ga$ be a connected algebraic group over a number field~$K$, and let
$S$ be a finite set of places of~$K$, containing the infinite places. If\/
$\Ga_\theints$ is Zariski dense in\/~$\Ga$, and the $S$-rank of every
$K$-simple factor of\/~$\Ga/ \Rad \Ga$ is at least two, then\/ $\Ga_\theints$
is a superrigid lattice in\/~$\GOS$.
 \qed
 \end{cor}

\begin{rem} \label{need-weak}
 Although any connected, noncompact simple Lie group~$G$ has no
finite-dimensional unitary representations, some of the lattices in~$G$
may have finite-dimensional unitary representations. Any such
lattice~$\Lambda$ is not a superrigid subgroup of~$G$, in the sense of
Defn.~\ref{super-def}. However, the Margulis Superrigidity Theorem
\cite[Thm.~IX.6.16(b), p.~332]{MargBook} asserts that $\Lambda$ satisfies
a weaker definition of superrigidity that considers only representations
$\alpha\colon \Lambda \to \GL_n(\real)$, such that
$\closure{\Lambda^\alpha}$ has no nontrivial connected, compact,
semisimple, normal subgroups.  If $\Lambda$ is solvable, then
$\closure{\Lambda^\alpha}$ is also solvable, so the assumption
on~$\closure{\Lambda^\alpha}$ is automatically satisfied and therefore is
irrelevant. However, the difference between the two definitions of
superrigidity {\it is} relevant in the context of Cor.~\ref{main-cor}.
Fortunately, the corollary is valid for either of the two definitions of
superrigidity.
 \end{rem}

\begin{ack}
 I would like to thank the mathematics departments of Hebrew University
and the Tata Institute of Fundamental Research for their hospitality; much
of this work was carried out during productive and enjoyable visits to
Jerusalem and Bombay. I am pleased to acknowledge very helpful
conversations with G.~A.~Margulis, Amos Nevo and Robert J.~Zimmer that
suggested my previous work should be generalized in this direction, and
pointed me in the right direction for a proof. I am also grateful to
J.--T.~Chang and Nimish Shah for discussions that clarified some of the
details.  This research was partially supported by a grant from the
National Science Foundation.
 \end{ack}

\section{An instructive proof of a special case}
 \label{easy-pf}

In this section, we present a simple proof of Thm.~\ref{my-thm}. This
special case illustrates many of the ideas involved in the proof of
our main theorem. In particular, this case illustrates the importance
of the existence of a syndetic hull.

\begin{defn}[({\cite[\S5]{Witte-super},
cf.~\cite[p.~6]{Fried-Goldman}})]
 \label{synd-def}
 Let $\Gamma$ be a closed subgroup of a Lie group~$G$. A {\it syndetic hull}
of~$\Gamma$ is a connected subgroup~$B$ that contains~$\Gamma$, such that
$B/\Gamma$ is compact.
 \end{defn}

Although our goal is a proof of Thm.~\ref{my-thm}, which deals only
with algebraic groups, it is convenient to prove a more general result
that applies to more general Lie groups \see{my-thm-Lie}, because it
is easier to work with simply connected groups, but the universal
cover of~$\Ga_{\real}$ is often \emph{not} an algebraic group.

The proof we give here is based on the same ideas
as~\cite{Witte-super}. However, the present proof is \emph{much} less
complicated because we do not bother to keep track of exactly when
it is necessary to pass to a finite-index subgroup or mod out a
finite group. In particular, we thereby avoid the need to discuss
nilshadows, which play an important role in~\cite{Witte-super}.

\begin{thm} \label{my-thm-Lie}
 Let\/ $\Gamma_1$ be a lattice in a solvable Lie group~$G_1$, and assume
$G_1$ has only finitely many connected components. If\/
$\Zi1{\Gamma_1} = \Gi1$, then\/ $\Gamma_1$ is superrigid in~$G_1$.
 \end{thm}

\begin{proof}
 Replacing $G_1$~and~$\Gamma_1$ by finite-index subgroups, we may
assume $G_1$ is connected. Then, by passing to the universal cover, we
may assume $G_1$ is simply connected.
 Given a homomorphism $\alpha \colon \Gamma_1 \to \GL_n(\real)$, let
$G_2 = \closure{\Gamma_1^\alpha}$. Replacing $\Gamma_1$ by a
finite-index subgroup, we may assume $G_2$ is connected, so
$[G_1,G_1]$ and $[G_2,G_2]$ are simply connected, nilpotent Lie
groups.

Let $G = G_1 \times G_2$ and $\Gamma = \graph(\alpha)$, so $\Gamma$ is
a discrete subgroup of~$G$. Any maximal compact torus of $\G$ is of
the form $T_1 \times S_1$, where $T_1$ and~$S_1$ are maximal compact
tori of $\Gi1$ and $\Ad G_2$, respectively. We may assume $T_1 \times
S_1$ contains a maximal compact torus~$T$ of $\Z{\Gamma}$. Because
$\Zi1{\Gamma_1} = \Gi1$, we know that the projection of~$T$ into
$\Gi1$ is all of~$T_1$. Therefore, $T S_1 = T_1 \times S_1$ is a
maximal compact torus of $\G$. Because $G_2$ is a real algebraic
group, there is a compact torus~$S$ of~$G_2$ with $\Ad_{G_2} S =
S_1$. Hence, Lem.~\ref{easy-synd} implies that some finite-index
subgroup of~$\Gamma$ has a simply connected syndetic hull~$X$.

It is not difficult  (cf.~\cite[Step~5 of pf.~of Thm.~6.4,
p.~174]{Witte-super}) to see that $X G_2 = G_1 G_2$ and $X \cap G_2 = e$,
from which it follows that $X$ is the graph of a homomorphism $\sigma
\colon G_1 \to G_2$. Because $X$ contains a finite-index subgroup
of~$\Gamma$, we know that $\sigma$ agrees with~$\alpha$ on a finite-index
subgroup of~$\Gamma_1$.
 \end{proof}

By definition, syndetic hulls are connected, so the proof of their
existence requires some way to prove that a subgroup is connected. The
following proposition plays this role in the proof of
Lem.~\ref{easy-synd} below. The proposition is proved
in~\S\ref{misc-sect}, but the statement is copied here for ease of
reference.

\begin{prop}[{\see{inverse-alm-conn}}] \label{inv-alm-conn-copy}
 Let $G$ be a connected, solvable Lie group, let
  $\rho \colon G \to \GL_n(\real)$ be a continuous homomorphism, and
let~$A$ be an almost-Zariski closed subgroup of\/ $\GL_n(\real)$. If
there is a compact, abelian subgroup~$S$ of~$G$ and a compact
torus~$T$ of~$A$ such that $S^\rho T$ is abelian and contains a
maximal compact torus of~$\closure{G^\rho }$, then the inverse image
$\rho ^{-1} (A)$ has only finitely many connected components. \qed
 \end{prop}

\begin{prop}[{\cite[Prop.~2.5, p.~31]{Raghunathan}}] \label{synd-nilp}
 Let\/ $\Gamma $ be a closed subgroup of a simply connected, nilpotent Lie
group~$G$.  Then\/ $\Gamma $ has a unique syndetic hull in~$G$. \qed
 \end{prop}

\begin{lem}[{\cite[Cor.~5.18]{Witte-super}}] \label{synd-abel-sc}
  If\/ $\Gamma $ is a discrete subgroup of a connected, abelian Lie
group~$G$, then some finite-index subgroup of\/~$\Gamma$ has a simply
connected syndetic hull in~$G$.
 \end{lem}

\begin{proof} Let $\cover{G}$ be the universal cover of~$G$. Replacing
$\Gamma$ by a finite-index subgroup, we may assume $\Gamma$ is
torsion-free (hence free abelian), so there is a
subgroup~$\cover{\Gamma}$ of~$\cover{G}$ that maps isomorphically
onto~$\Gamma$ under the covering map $\cover{G} \to G$. Let
$\cover{X}$ be the unique syndetic hull of~$\cover{\Gamma}$
in~$\cover{G}$ \see{synd-nilp}, and let $X$ be the image
of~$\cover{X}$ in~$G$. Then $X$ is a simply connected syndetic hull
of~$\Gamma$, as desired.
 \end{proof}

\begin{prop}[{\cite[Cor.~3.14]{Witte-super}}] \label{conn-norm-closed}
 Let $A$ be a connected Lie subgroup of a connected Lie
group~$G$. Then the normalizer of~$A$ in $\G$ is almost Zariski
closed.  \qed
 \end{prop}

\begin{prop}[{\cite[Cor.~3.10]{Witte-super}}] \label{conn-cent-closed}
 If $A$ is a subgroup of a connected, solvable Lie group~$G$ such that\/
$[G,G]$ is simply connected, then the centralizer of~$A$ in\/ $\G$ is
almost Zariski closed. \qed
 \end{prop}

The following lemma, which establishes the existence of an appropriate
syndetic hull, is a crucial ingredient in the proof of
Thm.~\ref{my-thm-Lie}. A generalization of this lemma is used in the
proof of our main theorem \see{simply-synd-norm}.

\begin{lem} \label{easy-synd}
 Let\/ $\Gamma$ be a discrete subgroup of a connected, solvable Lie
group~$G$, such that\/ $[G,G]$ is simply connected.
   If there is a compact, abelian subgroup~$S$ of~$G$ and a
compact torus~$T$ of\/ $\Z{\Gamma}$, such that\/ $(\Ad_G S) T$ is a
maximal compact torus of\/~$\G$, then some finite-index subgroup
of\/~$\Gamma$ has a simply connected syndetic hull in~$G$.
 \end{lem}

\begin{proof} Let $H = \Ad_G^{-1}(\Z{\Gamma})$.
Prop.~\ref{inv-alm-conn-copy} implies that~$H$ has only finitely many
connected components, so $H^\circ$ contains a finite-index subgroup
of~$\Gamma$. Thus, there is no harm in replacing~$G$ with~$H^\circ$,
which means we may assume $\Z{\Gamma} = \G$. 

From Prop.~\ref{synd-nilp}, we know that $[\Gamma,\Gamma]$ has a
unique syndetic hull~$U$ in $[G,G]$. The uniqueness implies that
$\Gamma$ normalizes~$U$. Then, because $U$ is connected and
$\Z{\Gamma} = \G$, we conclude that all of~$G$ normalizes~$U$
\see{conn-norm-closed}.  Thus, there is no harm in modding out~$U$,
so we may assume
$[\Gamma,\Gamma] = e$, that is, $\Gamma$ is abelian. So $\Gamma$
centralizes~$\Gamma$. Because $\Z{\Gamma} = \G$, this implies all
of~$G$ centralizes~$\Gamma$ \see{conn-cent-closed}, so $\Ad_G
\Gamma$ is trivial. Because $\Z{\Gamma} = \G$, this means $\Ad G$ is trivial,
so $G$ is abelian. The desired conclusion now follows from
Lem.~\ref{synd-abel-sc}.
 \end{proof}

 \section{Proof of the Main Theorem}
 \label{main-pf-sect}

This section presents a proof of the main theorem. However, instead of
the theorem as stated in~\S\ref{introduction}, we prove a more general
version that applies to groups whose algebraic structure is similar
to that of~$\GOS$ \see{super-thm}. We are not
particularly interested in this generalization for its own sake;
rather, the intention is to clarify the main ideas of the proof by
separating out the crucial hypotheses. The following proposition shows
that Main Theorem~\ref{main-thm} is indeed a special case of
Thm.~\ref{super-thm}. 

\begin{prop} \label{gos-is-ok}
 Let\/ $\Ga$ be a connected, solvable algebraic group over a number
field~$K$, and let $S$ be a finite set of places of~$K$, containing all the
infinite places, such that\/ $\Ga_\theints$ is Zariski dense in\/~$\Ga$.
 Let 
 \begin{itemize}
 \item $\phi\colon \Ga_S \to \Ga_\infty$ be the projection with
kernel\/~$\Ga_{\textnormal{f}}$;
 \item $G = \GOS$ \see{GOS-def};
 \item $\gi = G^\phi$;
 \item $\gf = \Ga_{\textnormal{f}}$; and
 \item $\Lambda = \hbox{any finite-index subgroup of~$\Ga_\theints$}$.
 \end{itemize}
Then
 \begin{enumerate}
 \item \label{1} $\gi$ is a solvable Lie group such that\/
$[\gi',\gi']=[\gi',(\gi')^\circ]$ is simply connected and nilpotent, for
every finite-index subgroup~$\gi'$ of~$\gi$;
 \item \label{2} $\gf$ is a locally compact, totally disconnected, solvable
group, such that\/ $[\gf,\gf]$ has no infinite discrete subgroups;
 \item \label{3} $G$ is an open subgroup of the direct product $\gi \times
\gf$;
 \item \label{4} $\Lambda$ is a lattice in $G$;
 \item \label{5} $\closure{\Ad_{\gi} \Lambda^\phi}$ is a finite-index
subgroup of\/ $\closure{\Ad \gi}$; and
 \item \label{6} $[\Lambda,\Lambda]$ is cocompact in\/ $[G,G]$.
 \end{enumerate}
 \end{prop}

\begin{proof}
 \pref{1} The definition of~$\gi$ implies that it is a (not
necessarily closed) Lie subgroup of~$\Ga_\infty$; so
$\gi$ is a solvable Lie group.
 We may write $\Ga$ as a semidirect product $\Ga = \mathbb T
\semiprod \mathbb U$, where $\mathbb T$ is a torus and $\mathbb U$ is
the unipotent radical. The unipotent group~$\mathbb U$ has no
nontrivial characters, so we have $\mathbb U_S \subset \GS$; in
particular, $\mathbb U_\infty \subset \GS$. It is well known that
$\mathbb U_{\script O}$ is a lattice in~$\mathbb U_\infty$
\see{is-lattice}, so the Borel Density Theorem \see{BDT-unip} implies
that
 $$ \Ad_{\Ga_\infty} \mathbb U_\infty \subset
 \closure{ \Ad_{\Ga_\infty} \mathbb U_{\script O} } \subset
 \closure{ \Ad_{\Ga_\infty} \Ga_\theints }
 . $$
 Therefore, $\mathbb U_\infty \subset \GOS$, so $\mathbb U_\infty \subset
G_\infty$. Because $\mathbb U_\infty$ is connected \see{unip-conn}, this
implies $\mathbb U_\infty \subset \gi'$, so $\gi' = T_\infty \semiprod
\mathbb U_\infty$, where $T_\infty = \mathbb T_\infty \cap \gi'$. Then,
because $\mathbb T$ is abelian, this implies that $[\gi',\gi'] = [T_\infty,
\mathbb U_\infty] [\mathbb U_\infty, \mathbb U_\infty]$. Because $\mathbb
U_\infty$ is connected, this implies $[\gi',\gi']$ is connected.
Therefore, being a connected subgroup of~$\mathbb U_\infty$, the commutator
subgroup $[\gi',\gi']$ is simply connected \see{solv-sc}.

\pref{2} The fact that $[\gf,\gf]$ has no infinite discrete subgroups
follows from the observation that $[\Ga,\Ga] \subset \mathbb U$ is
unipotent.

\pref{3} Because the range of each nonarchimedean valuation is a
discrete set (and the group of $K$-characters is finitely generated),
there is an open set $F \subset \Ga_{\textnormal{f}}$ such that
 $ ||\chi(x_v) ||_v = 1$ for all $x \in F$, all $K$-characters~$\chi$,
and all $v \in S-S_\infty$; hence $F \subset \GS$. Since
 $\Ad_{\Ga_\infty} F \subset \Ad_{\Ga_\infty}  \Ga_{\textnormal{f}} =
e$, then the definition of $\GOS$ implies that $F \subset \GOS$. 

Let $H$ be the identity component of~$\gi$. Because the range of each
nonarchimedean valuation is a discrete set, and $H$ contained in the
identity component of $(\GS)^\phi$, we must have $\prod_{v \in
S_\infty} || \chi(x_v) ||_v = 1$, for every $x \in H$, so $H \subset
\GS$. Furthermore, because $H \subset \gi = (\GOS)^\phi$, we must have
 $$ \Ad_{\Ga_\infty} H \subset
  \Ad_{\Ga_\infty} \gi \subset \closure{\Ad_{\Ga_\infty} \Ga_\theints} .$$
 Therefore, $H \subset \GOS$.

So $H \times F$ is an open subset of $\gi \times \gf$ contained in $\GOS =
G$. This establishes~\pref{3}.

\pref{4} Because $\Ga_\theints$ is a lattice in~$\GS$
\see{is-lattice} and $G$ is a closed subgroup of~$\GS$ that
contains~$\Ga_\theints$, we see that $\Ga_\theints$ is a lattice
in~$G$, as desired.

\pref{5} The desired conclusion is immediate from the fact that $\GOS
= G$.

\pref{6}
 Let $\mathbb U = \mathbb U^1 \supset \mathbb U^2 \supset \cdots \supset
\mathbb U^n = e$ be the descending central series of~$\mathbb U$, and let
$\mathbb U^i_\Lambda = \Lambda \cap \mathbb U^i_S$, a finite-index subgroup
of~$\mathbb U^i_\theints$. From Lem.~\ref{cocpct-lem} below, we see that
$[\mathbb U_\Lambda, \mathbb U_\Lambda^{n-2}]$ is cocompact in $\mathbb
U_S^{n-1}$. Then, by modding out $\mathbb U^{n-1}$ and applying the lemma
again, we see that $[\mathbb U_\Lambda, \mathbb U_\Lambda^{n-3}]$ is
cocompact in $\mathbb U_S^{n-2}$. Continuing in this manner, we see that
$[\mathbb U_\Lambda, \mathbb U_\Lambda]$ is cocompact in $[\mathbb U, \mathbb
U]_S$. Hence, there is no harm in modding out $[\mathbb U, \mathbb U]$, so we
may assume $\mathbb U$ is abelian.
 Then $[\Ga, \mathbb U, \mathbb U] = e$, so we conclude from
Lem.~\ref{cocpct-lem} below that $[\Lambda, \mathbb U_\Lambda]$ is
cocompact in $[\Lambda, \mathbb U]_S$.
 Hence, there is no harm in modding out $[\Lambda, \mathbb U]$, which
means that we may assume $\Lambda$~centralizes~$\mathbb U$. Because
$\Lambda$ is Zariski dense in~$\Ga$, this implies that $\Ga$
centralizes~$\mathbb U$, so $\Ga = \mathbb T \times \mathbb U$ is
abelian, so \pref{6}~is trivially true.
 \end{proof}

\begin{rem} From the proof, we see that it would suffice to make the
weaker assumption that $\Ad_{\Ga} \Ga_{\theints}$ is Zariski dense in
$\Ad \Ga$ in place of the assumption that $\Ga_\theints$ is Zariski
dense in~$\Ga$.
 \end{rem}

\begin{lem} \label{cocpct-lem}
 Let\/ $\Ga$ be an algebraic group over a number field~$K$, let $S$ be a
finite set of places of~$K$ that contains all the infinite places, and
let $\Lambda$ be a finite-index subgroup of\/~$\Ga_\theints$. If\/\/
$\mathbb X$~and\/~$\mathbb V$ are connected $K$-subgroups of\/~$\Ga$, such
that\/ $\mathbb V$~is unipotent, $\mathbb X$~normalizes\/~$\mathbb V$, and\/
$[\mathbb X, \mathbb V, \mathbb V] = e$, then\/ $[\Lambda \cap \mathbb
X_S, \mathbb \Lambda \cap \mathbb V_S]$ is cocompact in\/~$[\Lambda \cap
\mathbb X_S, \mathbb V]_S$.
 \end{lem}

\begin{proof}
 From the ascending chain condition on Zariski closed, connected subgroups, we
know there is a finite subset $\{x_1, \ldots, x_m\}$ of~$\Lambda \cap \mathbb X_S$ such that $[\Lambda \cap \mathbb X_S, \mathbb V] = [x_1, \mathbb V] [x_2,
\mathbb V] \cdots [x_m, \mathbb V]$. 
 The image of an $S$-arithmetic subgroup under a $K$-epimorphism is an
$S$-arithmetic subgroup \cite[Thm.~5.9, p.~269]{Plat-Rap}, so we see
that $[x_i, \Lambda \cap \mathbb V_S]$ is a finite-index subgroup of
$[x_i,\mathbb V]_\theints$, which is cocompact in $[x_i, \mathbb V]_S$
\see{is-lattice}.
 \end{proof}

\begin{defn}[{\cite[\S I, p.~8]{MargBook}}]
 Let $\Gamma$~and~$\Lambda$ be subgroups of a group~$G$. We say that
$\Lambda$ \emph{commensurabilizes}~$\Gamma$ if $\Gamma
\cap (\lambda^{-1} \Gamma \lambda)$ is a finite-index subgroup of both
$\Gamma$~and~$\lambda^{-1}\Gamma\lambda$, for every $\lambda \in
\Lambda$.
 \end{defn}

\begin{lem} \label{simply-synd-norm}
 Suppose $G$ is a solvable Lie group, and\/ $[G^\circ,G]$ is simply
connected.
 Let\/ $\Gamma$ be a discrete subgroup of~$G^\circ$ and let $\Lambda$ be a
subgroup of~$G$ that contains\/~$\Gamma$ and commensurabilizes it. Suppose
some subgroup of\/~$\Gamma$ has a syndetic hull~$U$, such that 
 \begin{itemize}
 \item $U$ is simply connected;
 \item $U \subset \nil G$;
 \item $U$ contains\/ $[\Lambda,\Lambda]$; and
 \item $U$ is normalized by~$\Lambda$.  
 \end{itemize}
  If there is a compact, abelian subgroup~$S$ of~$G^\circ$ and a compact
torus~$T$ of\/ $\Z{\Lambda}$, such that\/ $(\Ad_G S) T$ is a maximal compact
torus of~$\G$, then some finite-index subgroup of\/~$\Gamma$ has a simply
connected syndetic hull in~$G$ that contains~$U$ and is normalized
by~$\Lambda$.
 \end{lem}

\begin{proof}
 Because $\Go$ is a normal subgroup of $\G$, we know that $(\Ad_G S)T$
contains a maximal compact torus of $\Go$
\see{cpct-normal}, so $\Ad_{G^\circ}^{-1} (\Z{\Lambda})$ has
only finitely many connected components \see{inverse-alm-conn}. In
other words, letting $H = \Ad_{G}^{-1} (\Z{\Lambda})$, we know that
$H^\circ$~is a finite-index subgroup of $H \cap G^\circ$, so $H^\circ$
contains a finite-index subgroup of~$\Gamma$. Also, because $\Ad_G U$
is unipotent \see{nilrad-unip} and $\Gamma$ contains a cocompact
subgroup of~$U$, we conclude from the Borel Density Theorem
\pref{BDT-unip} that $\Ad_G U \subset \Z{\Gamma} \subset \Z{\Lambda}$,
so $U \subset H$. Hence, there is no harm in replacing $G$ with~$H$,
so we may assume $\Z{\Lambda} = \G$.

Because $\Lambda$ normalizes~$U$, and $\Z{\Lambda} = \G$, we conclude that all
of~$G$ normalizes~$U$ \see{norm-closed}.  Then there is no harm in modding
out~$U$, so $[\Lambda,\Lambda] = e$; that is, $\Lambda$ is abelian.

Because $\Gamma \subset \Lambda$, this implies that $\Lambda$
centralizes~$\Gamma$, so all of~$G$ centralizes~$\Gamma$
\see{cent-closed}. In particular, $G^\circ$ centralizes~$\Gamma$,
which means $\Gamma \subset Z(G^\circ)$.  Furthermore, $\Ad G^\circ
\subset
\Z{\Lambda}$ is abelian, so $G^\circ$ is nilpotent. This implies that
$Z(G^\circ)$ is connected \see{Z(G)-conn}. Hence, some finite-index
subgroup of~$\Gamma$ has a simply connected syndetic hull~$X$ in
$Z(G^\circ)$ \see{synd-abel-sc}.  All that remains is to show that $X$
is normalized by~$\Lambda$.

Let
 \begin{itemize}
 \item $\gio$ be the universal cover of~$G^\circ$,
 \item $\cover{X}$ be the connected subgroup of~$\cover{G}$ that covers~$X$,
 \item $\cover{\Gamma}$ be the inverse image of~$\Gamma$ in~$\cover{X}$, and
 \item $Z$ be the kernel of the covering map $\cover{G} \to G$.
 \end{itemize}
 Because $\Lambda$ centralizes~$\Gamma$, we know $[\cover{\Gamma},
\Lambda] \subset Z$. On the other hand, because $[G^\circ,G]$ is simply
connected, we know that $Z \cap [\cover{G}^\circ,G] = e$. Therefore,
$[\cover{\Gamma}, \Lambda] = e$, so $\Lambda$
normalizes~$\cover{\Gamma}$. Now $\cover{G}^\circ$ is a simply connected,
nilpotent Lie group, so $\cover{X}$ is the unique syndetic hull
of~$\cover{\Gamma}$ in~$\cover{G}^\circ$ \see{synd-nilp}, so we conclude that
$\Lambda$ normalizes~$\cover{X}$. Hence, $\Lambda$ normalizes~$X$, as
desired.
 \end{proof}

\begin{thm} \label{super-thm}
 Suppose
 \begin{itemize}
 \item $\gi$ is a solvable Lie group such that\/
$[\gi',\gi']=[\gi',(\gi')^\circ]$ is simply connected and nilpotent, for
every finite-index subgroup~$\gi'$ of~$\gi$; and
 \item $\gf$ is a locally compact, totally disconnected, solvable group,
such that\/ $[\gf,\gf]$ has no infinite discrete subgroups. 
 \end{itemize}
 Let 
 \begin{itemize}
 \item $G$ be an open subgroup of the direct product $\gi \times \gf$, and
 \item $\phi\colon \gi \times \gf \to \gi$ be the projection with kernel~$\gf$,
 \end{itemize}
 and let $\Lambda$ be a lattice in $G$ such that
 \begin{itemize}
 \item $\closure{\Ad_{\gi} \Lambda^\phi}$ is a finite-index subgroup
of\/ $\closure{\Ad \gi}$, and
 \item $[\Lambda_1,\Lambda_1]$ is cocompact in\/ $[G,G]$, for every
finite-index subgroup~$\Lambda_1$ of~$\Lambda$.
 \end{itemize}
 Then $\Lambda$ is superrigid in~$G$.

 More precisely,
 if $\alpha \colon \Lambda \to \GL_{n}(\real)$ is a homomorphism such that\/
 $$\closure{ \{ (\Ad_{\gi} \lambda, \lambda^{\alpha}) \mid \lambda \in
\Lambda \}}$$
 is connected,
 then there is a finite subgroup~$F$ of $Z(\closure{\Lambda^\alpha})$ such that
the induced homomorphism $\alpha_F\colon \Lambda \to
\closure{\Lambda^\alpha}/F$ extends to a continuous homomorphism
 $\sigma \colon G_1 \to \closure{\Lambda^{\alpha}}/F$, for some finite-index
subgroup~$G_{1}$ of~$G$.
 \end{thm}

\begin{proof}
 Replacing $G$~and~$\gi$ by finite-index subgroups, we may assume
 $\closure{\Ad_{\gi} \Lambda^\phi} = \closure{\Ad \gi}$.
 Assume for simplicity that $\gio$ is simply connected. (In the situation
of Proposition~3.1, this may be achieved by passing to a universal cover.)
 Let $K$ be a compact open subgroup of~$\gf$ contained in~$G$.
 Let $\Gamma=\Lambda\cap(\gio K)$; note that $\Lambda$
commensurabilizes~$\Gamma$.
 We may assume, by replacing $K$ with a finite-index subgroup, that
$\Gamma \cap K = e$ \see{residually-finite}, so $\phi$ is faithful
on~$\Gamma$.

\step{1}{There is a unique homomorphism $\beta\colon [\gi,\gi] \to
[\closure{\Lambda^{\alpha}},\closure{\Lambda^{\alpha}}]$ such that $\phi\beta$
agrees with $\alpha$ on $\Gamma \cap [\Lambda,\Lambda]$.}
 Let $\Gamma_1 = \Gamma \cap [\Lambda,\Lambda]$. Because $\phi$ is faithful
on~$\Gamma$, the homomorphism $\alpha|_{\Gamma_1} \colon \Gamma_1 \to
[\closure{\Lambda^{\alpha}},\closure{\Lambda^{\alpha}}]$ induces a homomorphism
$\bar{\alpha}\colon \Gamma_1^\phi \to
[\closure{\Lambda^{\alpha}},\closure{\Lambda^{\alpha}}]$. Because
$\Gamma_1^\phi$ is a lattice in $[\gi,\gi]$ \see{latt-comm}, we know that
$\bar{\alpha}$ extends to a unique homomorphism $\beta\colon [\gi,\gi] \to
[\closure{\Lambda^{\alpha}},\closure{\Lambda^{\alpha}}]$ \see{super-nilp}.

\step{2}{$\phi\beta$ also agrees with~$\alpha$ on $[\Lambda,\Lambda]$.}
 Because $[\gf,\gf]$ has no infinite discrete subgroups, we see that for every
$\lambda\in [\Lambda,\Lambda]$, there is some $n\in \integer^{+}$ with
$\lambda^{n} \in \gi\times K$. Because $\gi/\gio$ is abelian, we also know
$\lambda \in \gio \times \gf$, so we conclude that $\lambda^n \in \gio K$;
therefore, $\lambda^{n} \in\Gamma$.  Therefore, $\phi\beta$ agrees
with~$\alpha$ on~$\lambda^{n}$. Because $n$th roots are unique in a unipotent
Lie group such as $[\closure{\Lambda^{\alpha}},\closure{\Lambda^{\alpha}}]$, we
conclude that $\phi\beta$ agrees with~$\alpha$ on~$\lambda$, as desired.

\step{3}{$\beta$ extends to a homomorphism $\rho\colon\gio \to
\closure{\Lambda^{\alpha}}$ such that $\phi\rho$ agrees with $\alpha$ on a
finite-index subgroup of~$\Gamma$, and we have
 $ g^{\lambda^\phi \rho} = g^{\rho \lambda^\alpha}$ for all $g \in \gio$ and
$\lambda \in \Lambda$.}
 Let 
 \begin{itemize}
 \item $\hat G = \gi \times \closure{\Lambda^{\alpha}}$;
 \item $\hat\Gamma = \{ (\gamma^\phi, \gamma^\alpha) \mid \gamma \in \Gamma
\}$;
 \item $\hat\Lambda = \{ (\lambda^\phi, \lambda^\alpha) \mid \lambda \in
\Lambda \}$; and
 \item $U = \{ (u,u^\beta) \mid u \in [\gi,\gi] \}$.
 \end{itemize}
 From Lem.~\ref{simply-synd-norm}, we see that some finite-index subgroup
of~$\hat\Gamma$ has a simply connected syndetic hull~$\hat X$ in~$\hat G$, such
that $\hat X$ contains~$U$ and is normalized by~$\hat\Lambda$. It is not
difficult to see that $\hat X$ is the graph of a homomorphism $\rho\colon \gio
\to \closure{\Lambda^{\alpha}}$ (see Step~5 of the pf.~of
\cite[Thm.~6.4]{Witte-super}). Because $\hat\Lambda$ normalizes~$\hat X$, we
have
 $ g^{\lambda^\phi \rho} = g^{\rho \lambda^\alpha}$.

\step{4}{Completion of the proof.}
 Let $L = \graph(\alpha)$, $\hat X = \graph(\rho)$, and $X = \hat X K
[\gf,\gf]$. (So $L$, $\hat X$, and~$X$ are subgroups of $G \times
\closure{\Lambda^\alpha}$.) From Step~3, we see that $L$
normalizes~$X$; hence $XL$ is a subgroup of $G \times
\closure{\Lambda^\alpha}$.  Let $F = XL \cap (e \times
\closure{\Lambda^\alpha})$ and $H = \gio K [\gf,\gf]$. It suffices to
show that $F$~is a finite, normal subgroup of
$\closure{\Lambda^\alpha}$, for then $XL/F$ is the graph of a
well-defined homomorphism
 $$\sigma \colon H \Lambda \to \closure{\Lambda^\alpha}/F , $$
 and $H \Lambda$ is a finite-index subgroup of~$G$, because it is
an open subgroup that contains the lattice~$\Lambda$.

Because $L$ normalizes $X$, $L$, and~$(e \times
\closure{\Lambda^\alpha})$, it is obvious that $\Lambda^\alpha$
normalizes~$F$, so we only need to show that $F$~is finite.
 Because $[\Lambda,\Lambda]$ is cocompact in $[G,G]$ and $\Gamma$ is a
lattice in $\gio K$, we know that $\Gamma [\Lambda,\Lambda]$ is a
finite-index subgroup of $\Lambda \cap H$. Then, from Steps~2 and~3,
we conclude that $\phi\rho$ agrees with~$\alpha$ on a finite-index
subgroup of $\Lambda \cap H$, so $X$ contains a finite-index
subgroup of $L \cap (H \times \closure{\Lambda^\alpha})$. On the
other hand, because $X \subset H \times \closure{\Lambda^\alpha}$, we
have
 $$F
 = XL \cap (e \times \closure{\Lambda^\alpha})
 = \Bigl( X \bigl( L \cap (H \times \closure{\Lambda^\alpha})
\bigr) \Bigr)
 \cap (e \times \closure{\Lambda^\alpha}). $$
 Therefore, $X \cap (e \times \closure{\Lambda^\alpha})$ contains a
finite-index subgroup of~$F$. Because $X \cap (e \times
\closure{\Lambda^\alpha}) = e$, this implies that $F$~is finite, as
desired.
 \end{proof}

\section{Application to non-solvable groups}
 \label{non-solv-sect}

In this section, we prove Cor.~\ref{main-cor}. As described in
Rem.~\ref{need-weak}, we prove two versions of this corollary
\see{double-cor}.

\begin{defn}[{(cf.~\ref{super-def})}]
 Let $\Lambda$ be a subgroup of a topological group~$G$. Let us say
that $\Lambda$ is \emph{weakly superrigid} in~$G$ if, for every
continuous homomorphism $\alpha \colon \Lambda \to \GL_{n}(\real)$,
\emph{such that $\closure{\Lambda^\alpha}$ has no nontrivial
connected, compact, semisimple, normal subgroups},
 there are
 \begin{itemize}
 \item a finite-index open subgroup~$\Lambda_1$ of~$\Lambda$,
 \item a finite-index open subgroup~$G_{1}$ of~$G$, containing~$\Lambda_1$, and
 \item a finite, normal subgroup~$F$ of
$\closure{\Lambda_1^\alpha}$, where $\closure{\Lambda_1^\alpha}$ is
the (almost-)Zariski closure of $\Lambda_1^\alpha$ in
$\GL_n(\real)$, 
 \end{itemize}
 such that the induced homomorphism $\alpha_F\colon \Lambda_1 \to
\closure{\Lambda_1^\alpha}/F$ extends to a continuous homomorphism
 $\sigma \colon G_1 \to \closure{\Lambda_1^{\alpha}}/F$.
  \end{defn}

\begin{cor}[{(cf.~\ref{main-cor})}] \label{double-cor}
 Let\/ $\Ga$ be a connected algebraic group over a number field~$K$, and
let $S$ be a finite set of places of~$K$, containing the infinite places.
If\/ $\Ga_\theints$ is Zariski dense in~$\Ga$,
 and the image of\/ $\Ga_\theints$ in\/ $\GoS{(\Ga/\Rad\Ga)}$ is superrigid
\textnormal{(}or weakly superrigid\textnormal{)}, then\/ $\Ga_\theints$ is a
superrigid lattice in\/~$\GOS$ \textnormal{(}or a weakly superrigid lattice
in\/~$\GOS$, respectively\textnormal{)}.
 \end{cor}

\begin{proof}[Proof \textnormal{(cf.~pf.~of
\cite[Thm.~9.9]{Witte-super})}.]
 Suppose $\alpha \colon \Ga_{\theints} \to \GL_{n}(\real)$ is a
homomorphism. (If the goal is to prove that $\GOS$ is \emph{weakly}
superrigid, assume $\closure{\Ga_{\theints}^\alpha}$ has no
nontrivial connected, compact, semisimple, normal subgroups.)
 Let $\Gamma$ be a finite-index subgroup of $\Ga_{\theints}$, and let
$H = \closure{\Gamma ^\alpha }$. Replacing $\Gamma $ by a finite-index
subgroup, we may assume $H$~is connected and that there is a Levi
subgroup~${\La}$ of~${\Ga}$ such that $\Gamma = ({\La_S} \cap \Gamma)
\bigl({(\Rad\Ga)_S} \cap \Gamma \bigr)$. Let $L$~and~$\therad$ be
finite-index subgroups of $\GoS{\La}$ and $\GoS{(\Rad\Ga)}$,
respectively. Let $\alpha _L = \alpha |_{L \cap \Gamma }$ and
$\alpha_{\therad} = \alpha |_{\therad \cap \Gamma }$, and let $L_H =
\closure{(L \cap\Gamma)^\alpha }$.  Because $L \cap \Gamma $ is a
superrigid lattice in~$L$ (or weakly superrigid lattice,
respectively), it must be the case that $L_H$ is semisimple, so $L_H$
is a Levi subgroup of~$H$, and (after passing to a finite-index
subgroup) there are a finite (or compact, respectively), normal
subgroup~$F$ of $C_{L_H}(L^{\beta _L})$, and a continuous
homomorphism $\beta _L \colon L \to L_H/F$, such that
 $$ \gamma ^{\beta _L} =\gamma ^\alpha F , \qquad \forall
\gamma \in L \cap \Gamma \; .$$
 From Thm.~\ref{main-thm}, we also know that $\alpha _{\therad}$
extends to a homomorphism 
 $$\beta _{\therad} \colon \therad \to
\closure{(\therad \cap\Gamma)^\alpha } = \Rad H.$$

The semisimple group~$\La$ has no nontrivial $K$-characters, so
$\Gs{\La} = \La_S$. Therefore, $\La_\finite \subset \GoS{\La}$, so $L$
contains a finite-index subgroup of~$\La_\finite$. Thus, we may assume
$L = L_\infty \times L_\finite$, where $L_\infty = \La_\infty \cap L$
and $L_\finite = \La_\finite \cap L$.  Because $L_\finite$ is totally
disconnected and has no open normal subgroups of infinite index, we
know that $L_\finite^{\beta_L}$ is finite; therefore, replacing
$L$ by a finite-index subgroup, we may assume $\beta_L$ is
trivial on~$L_\finite$.

\CASE{1}{$C_H(\Rad H)$ has no nontrivial, compact, solvable, normal
subgroups.}
 (Note that if $\GOS$ has no nontrivial connected, compact,
semisimple, normal subgroups, then this implies that $C_H(\Rad H)$
has no compact normal subgroups at all, solvable or not.)
 In this case, the extension~$\beta _{\therad}$ is unique
(cf.~\cite[Cor.~6.11]{Witte-super}), so $\graph(\alpha )$ normalizes
$\graph(\beta _{\therad})$. 

Write $\Rad\Ga = \Ta \semiprod \Ua$, where $\Ta$ is a torus and $\Ua$
is the unipotent radical. We may assume $R = T \semiprod U$, where
$T$~and~$U$ are finite-index subgroups of~$\GoS{\Ta}$ and~$\Ua_S$,
respectively. Because $\graph(\alpha )$ normalizes $\graph(\beta
_{\therad})$, and $\La$ centralizes~$\Ta$, we see that
$\graph(\alpha_L)$ centralizes $\graph(\beta _{\therad}|_{T})$, so
$L_H$ centralizes $T^{\beta_\therad}$.

Let $\phi \colon L \to L_\infty$ be the projection with kernel
$L_\finite$, let $A = \{ (\gamma^\phi, \gamma^{\alpha}) \mid \gamma
\in L \cap \Gamma \}$, and let $\closure{A}$ be the almost-Zariski
closure of~$A$ in $L_\infty \times L_H$ \see{real-Zar}.

Let $U_\infty = \Ua_\infty \cap U$ and $U_\finite = \Ua_\finite \cap
U$, and assume $U = U_\infty \times U_\finite$.  Because $C_H(\Rad H)$
has no nontrivial, compact, solvable, normal subgroups, it is not
difficult to show that $U_\finite^{\beta_{\therad}}$ must be trivial.
Then, because $\graph(\alpha )$ normalizes $\graph(\beta
_{\therad})$, we see that $A$ must normalize
$\graph(\beta_{\therad}|_{U_\infty})$.  Thus, $\closure{A}$
normalizes $\graph(\beta _{\therad}|_{U_\infty})$
\see{norm-closed}. Therefore, $\closure{A} \cap (e \times L_H)$
centralizes~$U_\infty^{\beta _{\therad}} = U^{\beta _{\therad}}$.
Because $\closure{A} \cap (e \times L_H) \subset L_H$ also
centralizes~$T^{\beta _{\therad}}$, we conclude that $\closure{A}
\cap (e \times L_H)$ centralizes
 $\closure{(TU)^{\beta _{\therad}}} = \closure{R^{\beta _{\therad}}}
=\Rad H$
 \see{cent-closed}. Then, since
 $$ \closure{A} \cap (e \times L_H)
  \subset \graph(\beta _L) \cap (e \times L_H)
 = e \times F \; ,$$
 and $C_{H}(\Rad H)$ has no nontrivial compact, (solvable) normal
subgroups, we conclude that $\closure{A} \cap (e \times L_H)$ is
trivial.
 This means that $\closure{A}$ is the graph of a well-defined
homomorphism $\beta'_L \colon L_\infty \to L_H$ and, from the
definition of~$A$, we see that $\phi \beta'_L$
extends~$\alpha_L$. Hence, we may assume $\beta_L = \beta'_L$ (and
$F$ is trivial).
  Then $\graph(\beta_L)$ normalizes both
$\graph(\beta_{\therad}|_{U_\infty})$ and
$\graph(\beta_{\therad}|_{U_\finite})$ (the latter because
$\beta_{\therad}|_{U_\finite}$ is trivial), and, being contained in $L
\times L_H$, centralizes $\graph(\beta_{\therad}|_T)$. Hence,
$\graph(\beta _L)$ normalizes $\graph(\beta _{\therad})$, so the
function
 $$ \beta \colon L \semiprod \therad \to H \colon (l,r)
\mapsto l^{\beta _L} r^{\beta _{\therad}}$$
 is a homomorphism. This completes the proof in this case.

\CASE{2}{The general case.}
 Let $C$ be the (unique) maximal compact, solvable, normal subgroup of
$C_H(\Rad H)$ (which we no longer assume to be trivial).  From
Case~1, we know there is a finite-index subgroup~$G$ of $\GOS$ and a
homomorphism $\bar\beta \colon G \to H/C$ such that $\bar\beta $
extends the homomorphism induced by~$\alpha $.  Now $C^\circ$ is a
compact torus \see{cpct-abel}, so the Levi decomposition implies that
there is a normal subgroup~$J$ of~$H$ such that $JC = H$ and $J \cap
C$ is finite.  There is no harm in modding out this finite
intersection, so we may assume $J \cap C$ is trivial.
 Then $H/C$ is naturally isomorphic to~$J$, so we can think of
$\bar\beta $ as a homomorphism from~$G$ to~$J$. From the definition
of~$\bar\beta $, we have $\gamma ^\alpha \in \gamma ^\beta C$, for all
$\gamma \in \Gamma $. Since $C$ is central in~$H$ \see{cpct-central},
this implies that there is a homomorphism $\sigma \colon \Gamma \to C$
such that $\gamma ^\alpha = \gamma ^\beta \gamma ^\sigma$, for all
$\gamma \in \Gamma $. Because $C$ is abelian, we know that $\sigma$
is trivial on $[\Gamma ,\Gamma ]$; in particular, $\sigma $ is
trivial on a finite-index subgroup of $L \cap\Gamma$ and a
finite-index subgroup of $[\Gamma ,\therad \cap \Gamma ]$. From the
proof of Prop.~\ref{gos-is-ok}\pref{6}, we see that $[\Gamma ,\therad
\cap \Gamma ]$ is a cocompact subgroup of $[G,\therad]$, so, replacing
$\Gamma$ by a subgroup of finite index, we may assume that $\sigma $
is trivial on $(L \cap \Gamma ) ([G,\therad] \cap \Gamma ) = [G,G]
\cap\Gamma $. Thus, $\sigma $ extends to a homomorphism $\tau \colon G
\to C$ (for example, this follows by applying the main theorem
\pref{main-thm} to the abelian group $\Ga/[\Ga,\Ga]$). Then the
homomorphism $g \mapsto g^\beta g^\tau$ extends~$\alpha$, as desired.
 \end{proof}

\begin{rem} \label{real-Zar}
 $\Ga_\infty$ is of the form ${\mathbb A}_{\complex} \times {\mathbb
B}_{\real}$, where $\mathbb A$~and~$\mathbb B$ are algebraic groups
defined over~$\complex$ and~$\real$, respectively.
 By restriction of scalars, the $\complex$-points of an
$n$-dimensional algebraic group defined over~$\complex$ can be viewed
as the $\real$-points of a $2n$-dimensional algebraic group defined
over~$\real$. Thus, we see that $\Ga_\infty$ can be viewed as the
$\real$-points of an algebraic group defined over~$\real$. Therefore,
in a natural way, $\Ga_\infty$ has a Zariski topology.
 \end{rem}

\section{Miscellaneous facts from Lie theory}
 \label{misc-sect}

 In this section, we collect for convenient reference a number of
facts that are used in~\S\ref{main-pf-sect}.

 All locally compact groups (including all Lie groups) in this paper are
assumed to be second countable.

\subsection{Connected subgroups}

\begin{lem}[{(cf.~\cite[Thm.~VIII.1.1, p.~107]{Hochschild-alg})}]
\label{unip-conn}
 Every almost-Zariski closed, unipotent subgroup of\/ $\GL_n(\real)$ is
connected and simply connected.
 \qed
 \end{lem}

\begin{lem}[{\cite[Thm.~XII.2.2, p.~137]{Hochschild-Lie}}]
 \label{solv-sc}
  Every connected subgroup of any simply connected, solvable Lie group~$G$
is closed and simply connected. \qed
 \end{lem}

\begin{lem}[{(cf.~\cite[Cor.~I.5.3.7, p.~47]{Bourbaki})}]
\label{nilrad-unip}
 If $A$ is a connected subgroup of a Lie group~$G$, then $A \subset
\nil G$ if and only if\/ $\Z{A}$ is unipotent. In particular, $\Z{\nil
G}$ is unipotent. \qed
 \end{lem}

\begin{lem}[{\cite[Thm.~XVI.1.1, p.~188]{Hochschild-Lie}}]
\label{Z(G)-conn}
 If $G$ is a connected, nilpotent Lie group, then $Z(G)$ is connected. \qed
 \end{lem}

\begin{lem}[{\cite[Lem.~3.20]{Witte-super}}] \label{unip-closed}
 Every connected, unipotent Lie subgroup of $\GL_n(\real)$ is Zariski closed.
\qed
 \end{lem}

\subsection{Compact subgroups}

\begin{prop}[{\cite[Satz~4]{Iwa2}, \cite[Thm.~XIII.1.3,
p.~144]{Hochschild-Lie}}]
 \label{cpct-abel} 
 If a Lie group~$G$ is compact, connected, and solvable, then $G$ is
abelian. Hence, $G \iso \torus^n$, for some~$n$. \qed
 \end{prop}

\begin{lem}[{(cf.~\cite[Satz~5]{Iwa2})}]
 \label{cpct-central} 
 If $G$ is a connected Lie group, then every compact subgroup of\/
$\nil G$ is central in~$G$. \qed 
 \end{lem}

\begin{prop}[{\cite[Thm.~XV.3.1, pp.~180--181, and see
p.~186]{Hochschild-Lie}}] \label{cpct-conj}
 If $G$ is a Lie group that has only finitely many connected
components, then $G$ has a maximal compact subgroup~$K$, and every
compact subgroup of~$G$ is contained in a conjugate of~$K$. \qed
 \end{prop}

\begin{cor} \label{cpct-normal}
 Let $G$ be a Lie group that has only finitely many connected
components. If $H$ is a closed, normal subgroup of~$G$, and $K$ is a
maximal compact subgroup of~$G$, then $H \cap K$ is a maximal compact
subgroup of~$H$. \qed
 \end{cor}

\begin{cor} \label{tori-conj}
 Let $G$ be a connected, solvable Lie group. Then all the maximal compact
tori of\/ $\G$ are conjugate under\/~$\Ad G$.
 \end{cor}

\begin{proof}
 Write $\G = T \semiprod U$, where $T$ is a maximal torus of $\G$, and $U$
is the unipotent radical. Let $M = [T(\Ad G)] \cap U$. Then $\Ad G \subset
TM$, and $TM$ is almost Zariski closed \see{unip-closed}, so we must have
$TM = \G$. Thus, $M=U$, so $T (\Ad G) = \G$. All the maximal compact tori of
$\G$ are conjugate under $\G$  \see{cpct-conj}, so, because $T$ normalizes
(indeed, centralizes) the maximal compact torus that it contains, this
implies that the maximal compact tori are conjugate under $\Ad G$.
 \end{proof}

\begin{lem}[{\cite[Thm.~7.6, p.~61]{Hewitt-Ross}}]
\label{residually-finite}
 If $K$ is a totally disconnected, compact group, then $K$ is
residually finite.  \qed
 \end{lem}

\subsection{Representations of nilpotent groups}

\begin{thm}[{\cite[Thm.~2.11, p.~33]{Raghunathan}}] \label{super-nilp}
 Let $G$~and~$H$ be simply connected, nilpotent Lie groups, and let\/
$\Gamma$ be a lattice in~$G$. Then every homomorphism from\/~$\Gamma$
to~$H$ extends to a unique continuous homomorphism from~$G$ to~$H$. \qed
 \end{thm}

\begin{cor} \label{nilp-cent}
 Let $G$ be a simply connected, nilpotent Lie group. Then the trivial
automorphism is the only automorphism of~$G$ that centralizes a cocompact
subgroup of~$G$. \qed
 \end{cor}

The following is a useful special case of the Borel Density Theorem.

\begin{prop}[{\cite[Thm.~3.2, p.~45]{Raghunathan}}] \label{BDT-unip}
 If\/ $\Gamma$ is a lattice in a Lie group~$G$, and $\phi\colon G \to
\GL_n(\real)$ is a representation such that $G^\phi$ is unipotent, then\/
$\closure{\Gamma^\phi} = \closure{G^\phi}$. \qed
 \end{prop}

\subsection{$S$-arithmetic groups}

\begin{thm}[{\cite[Thm.~5.6, p.~264]{Plat-Rap}}] \label{is-lattice}
 Let\/ $\Ga$ be a connected, solvable algebraic group over a number
field~$K$, and let $S$ be a finite set of places of~$K$, containing
all the infinite places. Then\/ $\Ga_{\theints}$ is a lattice in\/~$\GS$
\seemore{Defn.}{GOS-def}. \qed
 \end{thm}

\begin{eg} \label{GOS-proper}
 $\GOS$ may be a proper subgroup of $\GS$, even if $\Ga_{\theints}$ is
Zariski dense in~$\Ga$. To see this, it suffices to construct an
anisotropic torus~$\Ta$ over a number field~$K$, such that
$\Ta_{\theints}$ is not Zariski dense in~$\Ta_\infty$
\see{real-Zar}. For then the desired example is obtained by forming a
semidirect product $\Ga = \Ta \semiprod \mathbb U$, such that
$C_{\Ta}(\mathbb U)$ is trivial, and letting $S = S_\infty$.

Let $p$~and~$q$ be two distinct primes, with $q \equiv 3 \pmod{4}$, let $K =
\rational(i,\sqrt{p})$, and let $\Ta = \SO(x^2 + q y^2)$. Then $\Ta$ is
defined over~$K$ (in fact, it is defined over~$\rational$) and is
$K$-anisotropic (because $q$~is not a square in~$K$).  Now $K$ is Galois
over~$\rational$ and has two places, both complex, so
 $$ \Ta_\infty = \SO(x^2 + q y^2)_{\complex} \times \SO(x^2 + q y^2)_
{\complex} .$$
 Let $\sigma$ be the nontrivial Galois automorphism of~$K$ that
fixes~$i$,  and let $\tau$ denote complex conjugation. Define $\phi \colon
\Ta_\infty \to \Ta_{\complex}$ by $(u,v)^\phi = u u^\tau v v^\tau$. Then
we have
 $$ (\Ta_{\script O})^\phi
 = \{(u, u^\sigma)^\phi \mid u \in \Ta_{\script O} \}
 = \{ u u^\tau u^\sigma u^{\sigma\tau} \mid u \in \Ta_{\script O} \}
 .$$
 Clearly, then, each element of $(\Ta_{\script O})^\phi$ belongs to
$\Ta_{\script O}$ and is fixed by the Galois group of~$K$. Therefore,
$(\Ta_{\script O})^\phi \subset \Ta_{\integer}$. But $\Ta_{\integer}$
is finite, because it is a discrete subset of the compact
group~$\Ta_{\real}$, so we conclude that $(\Ta_{\script O})^\phi$ is
not Zariski dense in~$\Ta$. Because $\phi$ is a surjective morphism of
algebraic groups (where $\Ta_\infty$ is endowed with the Zariski topology
described in Rem.~\ref{real-Zar}), this implies that $\Ta_{\script O}$
is not Zariski dense in~$\Ta_\infty$.
 \end{eg}

 \subsection{Centralizers and normalizers are almost Zariski closed}

The following two propositions are proved in \cite{Witte-super} only in the
case where $G$ is connected. The general case follows by applying essentially
the same proofs to the connected group $G^\circ \prodsemi \G^\circ$.

\begin{prop}[{\cite[Cor.~3.10]{Witte-super}}] \label{cent-closed}
 Let $A$ be a subgroup of~$G^\circ$, where $G$ is a solvable Lie group such
that\/ $[G^\circ,G]$ is simply connected. Then the centralizer of~$A$ in\/
$\G$ is almost Zariski closed. \qed
 \end{prop}

\begin{prop}[{\cite[Cor.~3.14]{Witte-super}}] \label{norm-closed}
 Let $A$ be a connected Lie subgroup of a Lie group~$G$. Then the normalizer
of~$A$ in\/ $\G$ is almost Zariski closed. In particular, $\Z{A}$
normalizes~$A$. \qed
 \end{prop}

\subsection{Virtual connectivity of an inverse image}

\begin{prop}[{\cite[Lem.~5.6]{Witte-super}}] \label{inverse-conn}
 Suppose $G$ is a connected, solvable Lie group, and $A$ is an
almost-Zariski closed subgroup of $\GL_n(\real)$. If $\rho \colon G \to
\GL_n(\real)$ is a continuous homomorphism, such that $A$ contains a maximal
compact torus of~$\closure{G^\rho }$, then the inverse image $\rho ^{-1}
(A)$ is a connected subgroup of~$G$. \qed
 \end{prop}

\begin{cor} \label{inverse-alm-conn}
 Let $G$ be a connected, solvable Lie group, let
  $\rho \colon G \to \GL_n(\real)$ be a continuous homomorphism, and let~$A$ be
an almost-Zariski closed subgroup of $\GL_n(\real)$. If there is a compact,
abelian subgroup~$S$ of~$G$ and a compact torus~$T$ of~$A$, such that
$S^\rho T$ is abelian and contains a maximal compact torus
of~$\closure{G^\rho }$, then the inverse image $\rho ^{-1} (A)$ has only
finitely many connected components.
 \end{cor}

\begin{proof} 
 By replacing $A$ with $A \cap \closure{G^\rho}$, and $T$ with $(T \cap
\closure{G^\rho})^\circ$, we may assume $A \subset \closure{G^\rho}$. By
replacing $T$ with a larger torus (and replacing $S$ with a conjugate
that commutes with this larger torus \see{tori-conj}), we may assume $T$ is a
maximal compact torus of~$A$.
 Let $H$ be the (unique) almost-Zariski closed, connected subgroup of
$\closure{G^\rho}$ that has $T$ as a maximal compact torus and satisfies
$H S^\rho = \closure{G^\rho}$; let $S_1 = S^\rho /(H \cap S^\rho)$. Now
$H$ contains the commutator subgroup of $\closure{G^\rho}$, so there is a
natural homomorphism $\closure{G^\rho} \to \closure{G^\rho}/H \iso S_1$.
Let $\sigma\colon G \to S_1$ be the composition of~$\rho$ with this
homomorphism.

Let $K$ be the kernel of~$\sigma$; we begin by showing that $K$ has only
finitely many connected components.  Replacing $S$ by a subgroup, we may assume
$K \cap S$ is finite. The fibration $K \cap S \to S \stackarrow{\sigma} S_1$
yields the following long exact sequence of homotopy groups \cite[Cor.~IV.8.6,
p.~187]{White}:
 $$ \pi_1(S) \to \pi_1(S_1) \to \pi_0( K \cap S ) . $$
 Because $K \cap S$ is finite, we conclude that the cokernel of the map 
 $\pi_1(S) \to \pi_1(S_1)$ is finite. Because $S \subset G$, this implies that
the cokernel of the map $\pi_1(G) \stackarrow{\sigma_*} \pi_1(S_1)$ is
finite. Thus, from the long exact sequence 
 $$ \pi_1(G) \stackarrow{\sigma_*} \pi_1(S_1) \to \pi_0(K) \to \pi_0(G) = 0 ,$$
 obtained from the fibration $K \to G \to S_1$, we conclude that $\pi_0(K)$ is
finite, as desired.

 Because $A \subset H$, it is easy to see that $\rho^{-1}(A) \subset K$. Thus,
the conclusion of the preceding paragraph implies that $K^\circ$ contains a
finite-index subgroup of $\rho^{-1}(A)$, so there is no harm in replacing~$G$
with~$K$, so we may assume $G^\sigma = e$, which means $S^\rho \subset H$, so
$A$ contains a maximal compact torus of~$\closure{G^\rho}$. Hence, the
proposition implies that $\rho^{-1}(A)$ is connected.
 \end{proof}

\begin{rem} $\rho^{-1}(A)$ need not be connected, even if the restriction
of~$\rho$ to~$T$ is faithful. For example, in $\torus^3$,
let
 \begin{itemize}
 \item $G = \{(e^{it}, e^{i\lambda t})\} \times \torus$, (where $\lambda$ is
irrational),
 \item $A = \torus \times \{(e^{2is}, e^{is})\}$, 
 \item $T = \{(1,1)\} \times\torus$,
 and
 \item $\rho =$ the inclusion $G \hookrightarrow \torus^3$.
 \end{itemize}
  Then $\rho^{-1}(A) = G \cap A$ has two components:
 $$ \{(e^{it}, e^{i\lambda t}, e^{i\lambda t/2})\}
 \qquad \textnormal{and} \qquad
 \{(e^{it}, e^{i\lambda t}, -e^{i\lambda t/2})\} . $$
 \end{rem}

\subsection{The commutator subgroup of a lattice}

\begin{lem} \label{latt-comm}
 Suppose
 \begin{itemize}
 \item  $G$ is a solvable Lie group such that $G^\circ$ is simply
connected, and\/ $[G,G] = [G,G^\circ]$ is nilpotent;
 \item $\Gamma $ is a subgroup of an abstract group~$\Lambda $, such that
$\Lambda $ commensurabilizes\/~$\Gamma $; and
 \item $\phi  \colon \Lambda \to G$ is a homomorphism such that\/
  $\Z{\Lambda^\phi } = \G$,
  the restriction of~$\phi $ to\/~$\Gamma $ is faithful, and\/
 $\Gamma ^\phi $ is a lattice in $G^\circ$.
 \end{itemize}
 Then\/ $(\Gamma\cap[\Lambda ,\Lambda ])^\phi $ is a lattice in\/
$[G,G]$.
 \end{lem}

\begin{proof}
 Let $U$ be the unique syndetic hull of
$(\Gamma\cap[\Lambda,\Lambda])^\phi$ in $[G,G]$ \see{synd-nilp}.
 Since $\Lambda$ commensurabilizes $\Gamma\cap[\Lambda ,\Lambda]$, the
uniqueness of~$U$ implies that $\Lambda^\phi$ normalizes~$U$. Because
$\Z{\Lambda^\phi }=\G$, this implies that $U$ is a normal subgroup of~$G$
\see{norm-closed}. There is no harm in modding out~$U$, so we may assume
$\Gamma\cap[\Lambda,\Lambda] = e$.

Let us now show that 
 \begin{equation} \label{fin-ind-cent}
 \textnormal{each $\lambda\in\Lambda$ centralizes a finite-index subgroup
of~$\Gamma $.} \tag{$*$}
 \end{equation}
 Because $\Lambda$ commensurabilizes~$\Gamma$, there is some finite-index
subgroup~$N$ of~$\Gamma$ with $N^{\lambda}\subset \Gamma$. Then
 $[N,\lambda] \subset \Gamma\cap[\Lambda,\Lambda] = e$,
 so $\lambda$ centralizes~$N$. 

By replacing~$\Gamma $ with a subgroup of finite index, we may assume that
$\Z{\Gamma ^\phi }$ is connected, in which case,
\pref{fin-ind-cent}~implies that $\Ad_G\Lambda^\phi $ centralizes
$\Z{\Gamma ^\phi }$. Because $\Z{\Lambda^\phi } = \G$, this implies that
$\Z{\Gamma ^\phi } \subset Z(\G)$. Therefore, because $G^\circ/\Gamma
^\phi $ is compact, we see that the image of~$G^\circ$ in $\G/Z(\G)$ is
compact. But compact, connected Lie groups are abelian \see{cpct-abel}, so
we conclude that $G^\circ$ is nilpotent. 

Therefore, \pref{fin-ind-cent}~implies that $\Lambda$ centralizes~$G^\circ$
\see{nilp-cent}. Then, because  $\Z{\Lambda^\phi } = \G$, we see that $G$
centralizes~$G^\circ$ \see{cent-closed}. Hence $[G,G] = [G,G^\circ] = e$,
so the desired conclusion is trivially true.
 \end{proof}

\section{Errata to \protect\cite{Witte-super}}

 \begin{itemize}
 \item In the second sentence of the abstract, $\Gamma$ should be assumed
to be discrete.
 \item The reference for Lem.~3.21 should be to \cite{Hochschild-alg},
which was mistakenly omitted from the bibliography.
 \item The references for Prop.~5.2 and Cor.~5.3 should be [R, Thm.~2.1,
p.~29] and [R, Prop.~2.5, p.~31], respectively.
 \item In the second paragraph of Step~5 of the proof of Thm.~6.4, the
reference should be to Prop.~5.4, not~5.2.
 \item In the statement of Prop.~6.10, one must assume $G_1^\beta \subset
\closure{\Gamma_1^\alpha}$.
 \item The proof of Prop.~6.10 should begin by noting that, because
$G_{1}^{\beta_{1}}$ is connected, there is no harm in assuming $G_{2}$ is
connected.
 \item There is an error in the proof of Thm.~9.9. Lines~4--6 of
page~191 (immediately following the displayed equation) should be
replaced with the following:
 ``and $C_{L_H}(\Rad H)$ has no compact, connected, normal subgroups, we
conclude that $F = \closure{\graph{\alpha}} \cap (e \times L_H)$ is
finite. Thus, $\bigl(\closure{\graph{\alpha}} \cap (L \times L_H)\bigr)/F$
is the graph of a well-defined homomorphism $\bar\beta'_L \colon L \to
L_H/F$. Because $L$ is algebraically simply connected, we can lift
$\bar\beta'_L$ to a homomorphism $\beta'_L\colon L \to L_H$. Note that
$\beta'_L$ agrees with~$\alpha$ on a finite-index subgroup of $L \cap
\Gamma$ (cf.~the argument in the final paragraph of this
proof). Therefore, by replacing $\beta_L$ with~$\beta'_L$, we may
assume $\sigma$ is trivial.
 In other words, $\beta _L$ extends $\alpha _L$.''
 \end{itemize}

\end{document}